\documentclass[11pt]{amsart}
\usepackage{amssymb}
\usepackage{tikz}
\newtheorem{lemma}{\bf Lemma}
\newtheorem{theorem}[lemma]{\bf Theorem}

\newtheorem{proposition}[lemma]{\bf Proposition}

\begin{document}
\parskip = 0mm
\title[A Question of Bj\"orner from 1976]{A ``Challenging Question'' of Bj\"orner from 1976: Every Infinite Geometric Lattice of Finite Rank Has a Matching}
\author{Jonathan David Farley}
\thanks{The author would like to thank Professor Anders Bj\"orner for sending him his two preprints.  Not only would the author {\sl like} to thank him, he {\sl does} thank him.}
\address{Department of Mathematics,
Morgan State University,
1700 E. Cold Spring Lane,
Baltimore, Maryland 21251, United States of America, {\tt lattice.theory@gmail.com}}



\begin{abstract} It is proven that every geometric lattice of finite rank greater than 1 has a matching between the points and hyperplanes.  This answers a question of P\'olya Prize-winner Anders Bj\"orner from the 1981 Banff Conference on Ordered Sets, which he raised as a ``challenging question'' in 1976.
\end{abstract}

\maketitle


\def\Qa{\mathbb{Q}_0}
\def\Qb{\mathbb{Q}_1}
\def\Q{\mathbb{Q}}
\def\card{{\rm card}}
\parskip = 2mm
\parindent = 10mm
\def\Part{{\rm Part}}
\def\P{{\mathcal P}}
\def\Eq{{\rm Eq}}
\def\cld{Cl_\tau(\Delta)}
\def\Csing{{\mathcal C}_{\{*\}}}
\def\Cftwo{{\mathcal C}_{{\rm fin}>1}}
\def\Cinf{{\mathcal C}_{\infty}}
\def\Pcf{{\mathcal P}_{\rm cf}}
\def\Fn{{\mathcal F}_n}
\def\proof{{\it Proof. }}


\vspace*{-4mm} 


At the famous 1981 Banff Conference on Ordered Sets---such luminaries as Erd\H{o}s, Professor Garrett Birkhoff, Dilworth, Turing Award-winner D. S. Scott, Daykin, A. Garsia, R. L. Graham, C. Greene, B. J\'onsson, E. C. Milner, and Oxford's H. A. Priestley attended---Bj\"orner asked (with MIT's Richard Stanley asking a question immediately afterwards, judging from the proceedings) if every geometric lattice $L$ of finite rank [$\geq2$] had a matching \cite[pp. xi, xii, and 799]{RivHB}.  Greene had proven this for finite lattices \cite[Corollary 3]{GreGJ}.  Bj\"orner had proven this in special cases \cite[Theorems 3 and 4]{BjoGF}---for modular lattices and for ``equicardinal lattices,'' i.e., lattices whose hyperplanes contained the same number of atoms.  In 1976, Bj\"orner wrote, ``It would be interesting to know if the result of our theorems 3 and 4 can be extended to all infinite geometric lattices, or at least to some classes of such lattices other than the modular and the equicardinal.'' In 1977, he proved it for lattices of rank 3 and for lattices of cardinality less than $\aleph_\omega$.  The P\'olya Prize-winner went on to ask at the Banff Conference if there exists a family $M$ of pairwise disjoint maximal chains in $L\setminus\{0,1\}$ whose union contains the set of atoms, saying, ``I showed this is true for modular $L$, and J. Mason showed it to be true for finite $L$.''  He conjectured this in 1977 (\cite[p. 18]{BjoGG}, \cite[p. 10]{BjoGF}), writing in 1976 that ``[a]nother challenging question, related to the existence of matchings, is whether maximal families of pairwise disjoint maximal proper chains do exist in infinite geometric lattices (cf. \cite{MasGC}).''

We answer Bj\"orner's 1976 question about matchings.

We selectively use some of the notation and terminology from \cite{DavPriJB} and \cite[Chapter II, \S8 and Chapter IV]{BirFG}.

Let $P$ be a poset. Let $x,y\in P$ be such that $x\leq y$.  The {\sl closed interval} $[x,y]$ is $\{z\in P: x\leq z\leq y\}$.  If $|[x,y]|=2$, we say $x$ is a {\sl lower cover} of $y$ and $y$ is an {\sl upper cover} of $x$ and denote it $x\lessdot y$.

Let $P$ be a poset with least element $0$.  An {\sl atom} or {\sl point} is a cover of $0$.  The set of atoms is $\mathcal A(P)$.  If $P$ is a poset with greatest element $1$, a {\sl co-atom}, {\sl co-point}, or {\sl hyperplane} is a lower cover of $1$.  The set of hyperplanes is $\mathcal H(P)$.

A poset is {\sl semimodular} if, for all $a,b,c\in P$, $a\lessdot b,c$ and $b\ne c$ imply there exists $d\in P$ such that $b,c\lessdot d$.  A {\sl geometric lattice of finite height} is a semimodular lattice $L$ with no infinite {\sl chains} (totally ordered subsets)---implying $L$ has a $0$ and a $1$---such that every element is a join of a subset of atoms.  It is known \cite[Theorem 9.4]{Nat} that such an $L$ is a complete lattice with a finite maximal chain and all maximal chains have the same size $r+1$, where $r$ is the {\sl height} or {\sl rank} of $L$. Moreover, every element is a join of a finite set of atoms and a meet of a subset of $\mathcal H(L)$ (see \cite[Lemma 1]{BjoGF}).  Every interval is a geometric lattice \cite[\S3.3, Lemma]{WelGF}.  The rank of $\downarrow x:=[0,x]$ is the {\sl rank} $r(x)$ of $x\in L$.  For $x,y\in L$, $r(x\vee y)+r(x\wedge y)\leq r(x)+r(y)$ \cite[Theorem 9.5]{Nat}.  For $x\in L$, let $\underline x:=\mathcal A(L)\cap\downarrow x$ and let $\overline x:=\mathcal H(L)\cap[x,1]$.

The following is a basic fact (see \cite[p. 3]{BjoGF}).

\begin{lemma} Let $L$ be a geometric lattice of finite height.  Let $a,b\in L$ be such that $a\leq b$.  Then any $x\in[a,b]$ has a {\rm modular complement} in $[a,b]$, i.e., there exists $y\in[a,b]$ such that $x\wedge y=a$, $x\vee y=b$, and $r(x)+r(y)=r(a)+r(b)$.
\end{lemma}

\proof If $x=c_0\lessdot c_1\lessdot\dots\lessdot c_k=b$, find $a_i\in\mathcal A(L)\cap\downarrow c_i\setminus \downarrow c_{i-1}$ for $i=1,\dots,k$.  Let $y=a\vee a_1\vee\dots\vee a_k$.  Clearly $r(y)-r(a)=k=r(b)-r(x)$, $x\vee y=b$, and $x\wedge y\geq a$.  As $r(a)\leq r(x\wedge y)\leq r(x)+r(y)-r(x\vee y)=r(a)+r(b)-r(b)=r(a)$, we have $x\wedge y=a$.\qed

See \cite[Chapters 2, 3, 5 and 8]{JecJF} and \cite[Appendix 2, \S3]{LanJB} for basic facts about ordinals and cardinals.  If $\kappa$ is a regular cardinal, a subset $\Omega\subseteq\kappa$ is {\sl closed in} $\kappa$ if for every non-empty subset $A\subseteq\Omega$, the supremum of $A$ is $\kappa$ or in $\Omega$; it is {\sl unbounded in} $\kappa$ if the supremum of $\Omega$ is $\kappa$; it is a {\sl club in} $\kappa$ if it is both.  A subset $\Omega\subseteq\kappa$ is {\sl stationary in} $\kappa$ if it intersects every club in $\kappa$; note that $|\Omega|=\kappa$.

We take our notation from \cite[\S\S2, 4, and 6]{AhaNasSheHC}. A {\sl society} is a triple $\Lambda=(M_\Lambda,W_\Lambda,K_\Lambda)$ where $M_\Lambda\cap W_\Lambda=\emptyset$ and $K_\Lambda\subseteq M_\Lambda\times W_\Lambda$.  If $A\subseteq M_\Lambda$ and $X\subseteq W_\Lambda$, then $K_\Lambda[A]:=\{w\in W_\Lambda: (a,w)\in K_\Lambda\text{ for some }a\in A\}$, and $\Lambda[A,X]:=\bigl(A,X,K_\Lambda\cap(A\times X)\bigr)$ is a {\sl subsociety} of $\Lambda$.  If $B\subseteq M_\Lambda$, then $\Lambda-B:=\Lambda[M_\Lambda\setminus B,W_\Lambda]$.  If $\Pi$ is a subsociety, then $\Lambda/\Pi:=\Lambda[M_\Lambda\setminus M_\Pi,W_\Lambda\setminus W_\Pi]$.  We call a subsociety $\Pi$ of $\Lambda$ {\sl saturated} if $K_\Lambda[M_\Pi]\subseteq W_\Pi$ and we denote this situation by $\Pi\lhd\Lambda$.

An {\sl espousal} for $\Lambda$ is an injective function $E:M_\Lambda\to W_\Lambda$ such that $E\subseteq K_\Lambda$.  A society is {\sl critical} if it has an espousal and every espousal is surjective.

If $I$ is a set and $\bar\Pi=(\Pi_i: i\in I)$ is a family of subsocieties of $\Lambda$, then $\bigcup\bar\Pi:=(\bigcup_{i\in I}M_{\Pi_i},\bigcup_{i\in I}W_{\Pi_i},\bigcup_{i\in I}K_{\Pi_i})$.
Assume $I$ is an ordinal.  If $\theta\leq I$, then $\bar\Pi_\theta$ denotes $(\Pi_i:i<\theta)$.  The sequence $\bar\Pi$ is {\sl non-descending} if $\Pi_i$ is a subsociety of $\Pi_j$ whenever $i<j<I$; it is {\sl continuous} if, in addition, $\bigcup\bar\Pi_\theta=\Pi_\theta$ for every limit ordinal $\theta<I$.  
If $I=J+1$, $\bar\Pi$ is a J-{\sl tower in} $\Lambda$ if $\bar\Pi$ is a continuous family of saturated subsocieties of $\Lambda$ such that $\Pi_0=(\emptyset,\emptyset,\emptyset)$.

Let $\Pi$ be a subsociety of $\Lambda$.  Assume $1\leq\kappa\leq\aleph_0$.  Then $\Pi$ is a $\kappa$-{\sl obstruction in} $\Lambda$ if $\Pi\lhd\Lambda$ and $\Pi-A$ is critical for some $A\subseteq M_\Pi$ such that $|A|=\kappa$.

Now assume $\kappa$ is a regular, uncountable cardinal.  A $\kappa$-tower $\bar\Sigma$ in $\Lambda$ is {\sl obstructive} if, for each $\alpha<\kappa$, $\Sigma_{\alpha+1}/\Sigma_\alpha$ is either (a) a $\mu$-obstruction in $\Lambda/\Sigma_\alpha$ for some $\mu<\kappa$ or (b) $(\emptyset,w,\emptyset)$ for some $w\in W_\Lambda$, and $\{\alpha<\kappa: \text{(a) holds at }\alpha\}$ is stationary in $\kappa.$ We say $\Pi$ is a $\kappa$-{\sl obstruction in} $\Lambda$ if $\Pi=\bigcup\bar\Sigma$ for an obstructive $\kappa$-tower $\bar\Sigma$ in $\Lambda$; by \cite[Lemmas 4.2 and 4.3]{AhaNasSheHC}, $\Pi\lhd\Lambda$.

For a society $\Lambda$, $\delta(\Lambda)$ is the minimum of $\{|B|: B\subseteq M_\Lambda\text{ such that }\Lambda-B\text{ has an espousal}\}$.

We will use the following theorems of Aharoni, Nash-Williams, and Shelah:

\begin{theorem}{\rm (from \cite[Lemma 4.2 and Corollary 4.9a]{AhaNasSheHC})} If $\Pi$ is a $\kappa$-obstruction, then $\delta(\Pi)=\kappa$.\qed
\end{theorem}

\begin{theorem}\cite[Theorem 5.1]{AhaNasSheHC} A society $\Lambda$ has an espousal if and only if it has no obstruction.\qed
\end{theorem}

We will say that a geometric lattice of finite rank $r\geq 3$ has a {\sl matching} if the society $\biggl(\mathcal A(L),\mathcal H(L),\le\cap\bigl(\mathcal A(L)\times\mathcal H(L)\bigr)\biggr)$ has an espousal.  (Since $\mathcal A(L)=\mathcal H(L)$ in geometric lattices of rank 2, we could say they also have a {\sl matching}.)

Greene proved:

\begin{theorem}\cite[Corollary 3]{GreGJ} Every finite geometric lattice of rank at least 2 has a matching.\qed
\end{theorem}

Bj\"orner proved:
\begin{theorem}\cite[Theorems 3 and 6]{BjoGG} Every geometric lattice of rank 3, or of finite height and cardinality less than $\aleph_\omega$, has a matching.\qed
\end{theorem}

We use the following results of Bj\"orner:

\begin{lemma} {\rm (\cite[Lemma 1]{BjoGG} and \cite[Theorem 1]{BjoGF})} Let $L$ be a geometric lattice of finite height.
(a) Let $p\in\mathcal A(L)$, $h\in\mathcal H(L)$ and assume $p\nleq h$.  Then $|\underline h|\le|\overline p|$.
(b) If $L$ is infinite, then $|\mathcal A(L)|=|\mathcal H(L)|=|L|$.\qed
\end{lemma}

\begin{theorem} \cite[Theorem 4]{BjoGG} Let $L$ be an infinite geometric lattice of finite height such that $|\downarrow x|<|L|$ for every $x\in L$ of rank 2. If $|L|$ is a regular cardinal, then $L$ has a matching.\qed
\end{theorem}

Bj\"orner also uses this theorem of Milner and Shelah:

\begin{theorem} \cite[Theorem]{TveGF} Let $\Gamma=(M,W,K)$ be a society such that $K[m]\ne\emptyset$ for all $m\in M$ and such that $(m,w)\in K$ implies $|K^{-1}[w]|\le |K[m]|$.  Then $\Gamma$ has an espousal.\qed
\end{theorem}

We are ready to begin answering Bj\"orner's question.

\begin{lemma} Let $L$ be a geometric lattice of finite height.  Let $B\subseteq\mathcal A(L)$.  Let $\mathcal L(B)$ be the subposet $\big\{\bigvee_L\{b_1,\dots,b_n\}:n\in\mathbb N_0,b_1,\dots,b_n\in B\big\}$.

Then $\mathcal L(B)$ is a geometric lattice of finite height with rank $r_L(\bigvee_L B)$, and $\mathcal A\big(\mathcal L(B)\big)=B$.  The inclusion map is order- and cover-preserving.  Also, $0_{\mathcal L(B)}=0_L$ and $|\mathcal L(B)|$ is either finite or $|B|$.  If $1_{\mathcal L(B)}=1_L$, then $\mathcal H\big(\mathcal L(B)\big)\subseteq\mathcal H(L)$.
\end{lemma}

\proof Since $\mathcal L(B)$ is closed under arbitrary joins, it is a complete lattice (e.g., \cite[Theorems 2.31 and 2.41]{DavPriJB}).  Letting $n$ equal $0$ or $1$, we get $\{0_L\}\cup B\subseteq\mathcal L(B)$ and so $B\subseteq\mathcal A\big(\mathcal L(B)\big)$.  But for $n\ge2$, $b_1\vee b_2\vee\dots\vee b_n\ge b_1$, so $\mathcal A\big(\mathcal L(B)\big)\subseteq B$.  Clearly every element of $\mathcal L(B)$ is a join of atoms.  Let $m,n\in\mathbb N_0$ and let $b_1,\dots,b_n,c_1,\dots,c_m\in B$.  Assume $b_1\vee\dots\vee b_n\lessdot_{\mathcal L(B)} c_1\vee\dots\vee c_m$.  Then $m\ge1$.  Pick $r\in\{1,\dots,m\}$ such that $b_1\vee\dots\vee b_n< b_1\vee\dots\vee b_n\vee c_r\in\mathcal L(B)$.  Then $ b_1\vee\dots\vee b_n< b_1\vee\dots\vee b_n\vee c_r\le b_1\vee\dots\vee b_n\vee c_1\vee\dots\vee c_r\vee\dots\vee c_m=c_1\vee\dots\vee c_m$.  As $ c_1\vee\dots\vee c_m$ covers $ b_1\vee\dots\vee b_n$ in $\mathcal L(B)$, we conclude $b_1\vee\dots\vee b_n\vee c_r=c_1\vee\dots\vee c_m$.  By semimodularity in $L$, $b_1\vee\dots\vee b_n\lessdot_L b_1\vee\dots\vee b_n\vee c_r=c_1\vee\dots\vee c_m$.  Now let $k\in\mathbb N_0$ and let $d_1,\dots,d_k\in B$.  Assume that $b_1\vee\dots\vee b_n\lessdot_{\mathcal L(B)} d_1\vee\dots\vee d_k$ and $c_1\vee\dots\vee c_m\ne d_1\vee\dots\vee d_k$.  As before, for some $s\in\{1,\dots,k\}$, $ b_1\vee\dots\vee b_n \lessdot_L b_1\vee\dots\vee b_n\vee d_s=d_1\vee\dots\vee d_k$.  Thus $c_r\nleq d_1\vee\dots\vee d_k$ and $d_s\nleq c_1\vee\dots\vee c_m$.  By semimodularity, $c_1\vee\dots\vee c_m\lessdot_L c_1\vee\dots\vee c_m\vee d_s=b_1\vee\dots\vee b_n\vee c_r\vee d_s=d_1\vee\dots\vee d_k\vee c_r$ and $ d_1\vee\dots\vee d_k\lessdot_L d_1\vee\dots\vee d_k\vee c_r$; hence $c_1\vee\dots\vee c_m$, $d_1\vee\dots\vee d_k\lessdot_{\mathcal L(B)} b_1\vee\dots\vee b_n\vee c_r\vee d_s$.  This shows that $\mathcal L(B)$ is a geometric lattice, of finite height since $L$ has no infinite chains, with $1_{\mathcal L(B)}=\bigvee_L B$.  As $\bigvee_L B=\bigvee_L\{b_1,\dots,b_n\}$ for some $n\in\mathbb N_0$ and some $b_1,\dots,b_n\in B$, picking the smallest such $n$ and using semimodularity in $L$ and $\mathcal L(B)$, we see that $r_L(\bigvee_L B)=r_{\mathcal L(B)}(\bigvee_L B)$, namely $n$.  If $1_{\mathcal L(B)}=1_L$, then the hyperplanes of $L$ and $\mathcal L(B)$ have the same rank; thus $\mathcal H\big(\mathcal L(B)\big)\subseteq\mathcal H(L)$.

The cardinality of $\mathcal L(B)$ follows from standard arguments (or see \cite[Theorem 1]{BjoGF}).\qed

\begin{proposition} Let $\lambda$ be a singular cardinal.  Assume that every geometric lattice of finite rank at least $2$ and of cardinality less than $\lambda$ has a matching.  Then every geometric lattice of finite rank at least $2$ of cardinality $\lambda$ has a matching.
\end{proposition}

\proof (Compare this with the proof of \cite[Theorem 6.4]{AhaNasSheHC}.)  Assume not, for a contradiction. Then by Theorem 3, the society $\Gamma=\bigg(\mathcal A(L),\mathcal H(L),\le\cap\big(\mathcal A(L)\times\mathcal H(L)\big)\bigg)$ has a $\kappa$-obstruction $\Pi=(M,W,K)$, where $L$ is the lattice (and $L$ has rank at least $3$).  Since $|M|\le\lambda$, then by Theorem 2, we have $\kappa\le\lambda$---indeed $\kappa<\lambda$, since $\kappa$ is finite or a regular cardinal.  By Theorem 2, there exists $A\subseteq M$ such that $|A|=\kappa$ and $\Pi-A$ has an espousal, $H$.  Let $R\subseteq\mathcal A(L)$ be a finite subset such that $1_L=\bigvee R$.

Let $B_0=A\cup R$, and, for $n<\omega$, if $B_n$ is defined, let $B_{n+1}=B_n\cup H^{-1}\bigg(\mathcal H\big(\mathcal L(B_n)\big)\bigg)$.  Note that $R\subseteq B_n$ for all $n<\omega$, so the rank of $\mathcal L(B_n)$ is the rank of $L$ and $\mathcal H\big(\mathcal L(B_n)\big)\subseteq\mathcal H(L)$ by Lemma 9.

Let $B=\bigcup_{n<\omega} B_n\subseteq M\cup R$.  Now $|B_0|\le\max\{\kappa,\aleph_0\}<\lambda$.  If $n<\omega$ and $|B_n|\le\max\{\kappa,\aleph_0\}$, then $|\mathcal H\big(\mathcal L(B_n)\big)|\le\max\{\kappa,\aleph_0\}$, so $|B_{n+1}|\le\max\{\kappa,\aleph_0\}+\max\{\kappa,\aleph_0\}=\max\{\kappa,\aleph_0\}$.  Hence $|B|\le\aleph_0\max\{\kappa,\aleph_0\}=\max\{\kappa,\aleph_0\}<\lambda$.

As $R\subseteq B$, Lemma 9 shows that $|\mathcal L(B)|<\lambda$ and $\mathcal H\big(\mathcal L(B)\big)\subseteq\mathcal H(L)$, so $\mathcal L(B)$ has a matching. Let $G$ be the espousal.  Since 
$$
H[(M\setminus A)\setminus(M\setminus A)\cap B]\cap\mathcal H\big(\mathcal L(B)\big)=\emptyset
$$
and $A\subseteq B\cap M$---so that $M=[(M\setminus A)\setminus(M\setminus A)\cap B]\cup(B\cap M)$---we know $H|_{(M\setminus A)\setminus(M\setminus A)\cap B}\cup G|_{B\cap M}$ is an espousal of $\Pi$, as $\Pi\lhd\Gamma$, contradicting Theorem 2.\qed

With Theorem 5, Proposition 10 extends Bj\"orner's work to $\aleph_\omega$.  But using the argument of \cite[Theorem 6]{BjoGG} almost verbatim, we can settle Bj\"orner's first question from the 1981 Banff Conference on Ordered Sets.  Bj\"orner already did the heavy lifting in proving Theorem 5, but to make it clear that his proof is what we need, we include it.

\begin{theorem} Every geometric lattice of finite rank greater than $1$ has a matching.
\end{theorem}

\proof The proof is drawn from \cite[pp. 10--13]{BjoGG}. Assume we have a counterexample $L$ of smallest cardinality, and, among those counterexamples, one of smallest rank.  By Theorems 4 and 5 and Proposition 10, we can assume $|L|$ is a regular cardinal and that $L$ has rank at least $4$.  By Theorem 7, there is $\ell_0\in L$ of rank $2$ such that $|\downarrow\ell_0|=|L|$.

Assume that $|\overline p|=|L|$ for all $p\in\underline{\ell_0}$.  Consider any $q\in\mathcal A(L)\setminus\underline{\ell_0}$ and consider the rank $3$ geometric lattice $\downarrow(q\vee\ell_0)$.  By Lemma 6(b), $$
|L|=|\downarrow{\ell_0}|=|\underline{\ell_0}|\le|\{c\in\downarrow(q\vee\ell_0): q\lessdot c\}|
$$
(by Lemma 6(a))$=|\overline q|$ (by Lemma 6(b)), so $|L|=|\overline q|$.  Hence $|\overline p|=|L|$ for all $p\in\mathcal A(L)$.  By Theorem 8, $L$ has a matching.

So now assume $|\overline q|<|L|$ for some $q\in\underline{\ell_0}$.

{\it {\bf Case 1.} Every cover of $q$ except $\ell_0$ covers only one other atom.}

Then define $s:\mathcal A(L)\setminus{\underline{\ell_0}}\to\{x\in L:q\lessdot x\}$ by $s(p)=p\vee q$ for all $p\in\mathcal A(L)\setminus{\underline{\ell_0}}$.  In this case, $s$ is one-to-one.  By the minimality of $L$, the geometric lattice $\uparrow q$ has a matching $t: \{x\in L:q\lessdot x\}\to\overline q$.  We will define a matching $f$ for $L$.

Let $f(p):=t\big(s(p)\big)$ for all $p\in\mathcal A(L)\setminus\underline{\ell_0}$ and let $f(q):=t(\ell_0)$; we just need to define $f$ on $\underline{\ell_0}\setminus\{q\}$.  Pick $h_0\in\overline{\ell_0}$ and let $z$ be a modular complement of $\ell_0$ in $\downarrow h_0$.  Define $R:\underline{\ell_0}\to\{x\in L:z\lessdot x\lessdot h_0\}$ by $R(p)=p\vee z$ for all $p\in\underline{\ell_0}$.  This function is one-to-one: If $p,p'\in\underline{\ell_0}$ but $p\ne p'$ and $p\vee z=p'\vee z$, then $p\vee z=p\vee p'\vee z=\ell_0\vee z=h_0$, a contradiction.

If $p\in\underline{\ell_0}\setminus\{q\}$, then $q\nleq R(p)$ (or else $R(p)=p\vee q\vee z=\ell_0\vee z=h_0$, a contradiction), so $R(p)$ is covered by exactly one hyperplane in $\overline q$, namely $q\vee R(p)$, and this is $h_0$.  Since $f[\big(\mathcal A(L)\setminus\underline{\ell_0}\big)\cup\{q\}]\subseteq\overline q$, if $p\in\underline{\ell_0}\setminus\{q\}$, we can let $f(p)$ be any hyperplane covering $R(p)$ except $h_0$.  If $p_1,p_2\in\underline{\ell_0}\setminus\{q\}$ but $p_1\ne p_2$ and $f(p_1)=f(p_2)$, then $f(p_1)$ covers $R(p_1)=p_1\vee z$ and covers $R(p_2)=p_2\vee z$, so $f(p_1)=p_1\vee p_2\vee z=\ell_0\vee z=h_0$, a contradiction.  Thus $f$ is a matching.

{\it {\bf Case 2.} There exists $\ell_1\in L\setminus\{\ell_0\}$ such that $q\lessdot\ell_1$ and $|\underline{\ell_1}|\ge3$.}

Let $p_1,p_2\in\underline{\ell_1}$ be such that $|\{p_1,p_2,q\}|=3$.  Since $q=\ell_0\wedge\ell_1$, we have $p_1,p_2\nleq\ell_0$.  Let $h_0\in\overline{\ell_0}$ be such that $p_1\nleq h_0$.  (Pick a modular complement of $\ell_0\vee p_1$ in $\uparrow\ell_0$.)  If $p_2\le h_0$, then $q\le\ell_0\le h_0$ implies $\ell_1=p_2\vee q\le h_0$, and so $p_1\le h_0$, a contradiction.  Hence $p_2\nleq h_0$.

By the minimality of $L$, $\downarrow h_0$ has a matching $g:\underline{h_0}\to C:=\{x\in L: x\lessdot h_0\}$.  Let $C_2:=\{c\in C: |\overline c|=2\}$ and let $C_3:=C\setminus C_2$.  We will show that $|C_3|=|L|$.

Because $\{p_1,p_2\}\in\mathcal A(L)\setminus\underline{h_0}$, we have that $q\le\ell_1=p_1\vee p_2\le\bigvee\mathcal A(L)\setminus\underline{h_0}=:y$.

{\it Claim.  For $x\in C$, $x\in C_2$ if and only if $x=h_0\wedge h$ for some $h\in\overline y$.}

{\it Proof of claim.} We have a partition of $\mathcal A(L)\setminus\underline x$: $\{\underline k\setminus\underline x: x\lessdot k\in L\}$.

If $x\in C_2$, then $x\lessdot h$ for some $h\in\mathcal H(L)\setminus\{h_0\}$ and so $x=h_0\wedge h$.  If $w\in\mathcal A(L)\setminus\underline{h_0}$, then $w\notin\underline x$, so $x\lessdot w\vee x\in\mathcal H(L)$ but $w\vee x\ne h_0$, so $w\vee x=h$.  Hence $w\le h$.  Therefore $y\le h$.

Conversely, if $x=h_0\wedge h$ for some $h\in\overline y$, then $h_0\ne h$. If there exists $h'\in\overline x\setminus\{h_0,h\}$, then, for some $a\in\underline{h'}\setminus\underline x$, $h'=a\vee x$.  Hence $a\notin \underline{h_0}\setminus\underline x$, and thus $a\le y$, so $a\le h$ and thus $a\le h\wedge h'=x$, a contradiction.  Hence $x\in C_2$.\qed

By the claim, $|C_2|\le|\overline y|$.  But $q\le y$ implies that $\overline y\subseteq\overline q$ and, since $|\overline q|<|L|$, we conclude $|C_2|<|L|$.  By Lemma 6(b), $|C|=|L|$, so $|C_3|=|L|$.

We now define our matching as follows: Since $|\mathcal A(L)\setminus\underline{h_0}|\le|C_3|$, take any injection $b:\mathcal A(L)\setminus\underline{h_0}\to C_3$ and let $f(p)=p\vee b(p)$ for $p\in\mathcal A(L)\setminus\underline{h_0}$.  For $p\in\underline{h_0}$, let $f(p)$ be any cover of $g(p)$ except $h_0$ or, in case $g(p)=b(p')$ for some $p'\in\mathcal A(L)\setminus\underline{h_0}$, except $f(p')$.  (We can do this since $b(p')\in C_3$.)  If $x',x''\in C$ and $x'\ne x''$, then $\uparrow x'\cap\uparrow x''=\uparrow h_0$; hence if $p,p'\in\mathcal A(L)\setminus\underline{h_0}$ and $p\ne p'$ but $f(p)=f(p')$ (so $p\vee b(p)=p'\vee b(p')$), then $f(p)\ge h_0$; but $r\big(f(p)\big)=r(h_0)$, so $f(p)=h_0$ and $p\le h_0$, a contradiction.

If for some $p\in\underline{h_0}$ and $p'\in\mathcal A(L)\setminus\underline{h_0}$ we have $f(p)=f(p')$, then $g(p)\lessdot p'\vee b(p')$.  Since $g(p),b(p')\lessdot h_0$, then $\{g(p),b(p'),h_0,p'\vee b(p')\}$ would be a 4-element crown (also called a ``cycle'') of elements in consecutive ranks---impossible in a lattice---unless $g(p)=b(p')$, which we have ruled out.

If $p,p'\in\underline{h_0}$ and $f(p)=f(p')$ but $p\ne p'$, then $g(p)\ne g(p')$ and $g(p),g(p')\lessdot h_0$ and $f(p)$ is a cover of $g(p),g(p')$ distinct from $h_0$, so we get another impossible 4-crown.

Hence $f$ is one-to-one.\qed

This answers the question of Bj\"orner from 1976 that was the first question he stated at the 1981 Banff Conference on Ordered Sets.  A good approach to the second would be to use \cite{AhaBerJI} and \cite{LogShaJD}; the latter contains a theorem that, when he first read it, made this writer feel that it could hold its own alongside many classic results in combinatorics.


\end{document}